\documentclass[letterpaper, 10 pt, conference]{ieeeconf}  

\IEEEoverridecommandlockouts                              
\usepackage[all]{xy}
\usepackage{color,url}
\usepackage{subfigure}
\usepackage{bm,amsmath}
\usepackage{marvosym}
\usepackage{graphicx}
\usepackage{wrapfig}
\usepackage{epsf}
\usepackage{times,mathptm}
\usepackage{url}
\usepackage{algorithmic}
\usepackage{algorithm2e}
\usepackage{graphicx}
\graphicspath{Figures/}
\newcommand{\pv}{\mathbf{p}}
\newcommand{\pvs}{\pv^s}

\newcommand{\One}{\mathbf{1}}

\newcommand{\sv}{\mathbf{s}}

\newcommand{\fv}{\mathbf{f}}

\newcommand{\Gr}{\mathcal{G}}
\newcommand{\So}{\mathbf{S}}
\newcommand{\V}{\mathcal{V}}
\newcommand{\E}{\mathcal{E}}
\newcommand{\N}{\mathrm{Neb}}

\newcommand{\ind}{\mathbf{x}}
\newcommand{\Lap}{\mathbf{L}}

\newcommand{\flim}{\bar{\fv}}
\newcommand{\pvlim}{\bar{\pv}}
\newcommand{\svlim}{\bar{\sv}}
\newcommand{\pvslim}{\bar{\pvs}}

\newcommand{\pvr}{\pv^r}
\newcommand{\pvg}{\pv^g}
\newcommand{\pvl}{\pv^l}

\newcommand{\popt}{{\pvs}^*}
\newcommand{\sopt}{{\sv}^*}
\newcommand{\sr}{{\sv}^r}
\newcommand{\abs}[1]{\left|{#1}\right|}

\newcommand{\cost}{\mathbf{c}}
\newcommand{\tran}[1]{{#1}^{\bf T}}
\newcommand{\costg}{\cost_g}
\newcommand{\costs}{\cost_s}
\newcommand{\costsp}{\cost_s^p}
\newcommand{\pvglim}{\bar{\pv_{gr}}}

\date{}

\title{\Large \bf Storage Sizing and Placement\\ through Operational and Uncertainty-Aware Simulations}

\author{
{\bf Krishnamurthy Dvijotham}$^{(1)}$, {\bf Scott Backhaus} $^{(2)}$ and {\bf Misha Chertkov} $^{(3)}$\\
Department of Computer Science and Engineering, University of Washington, Seattle, WA 98195 $^{(1)}$\\
Material Physics and Applications Division $^{(2)}$ and \\
Center for Nonlinear Studies and Theoretical Division, $^{(3)}$\\
LANL, Los Alamos, NM 87545, USA}

\begin{document}

\maketitle
\thispagestyle{empty}
\pagestyle{empty}

\begin{center}
{\bf \large Abstract}
\end{center}
{\it
As the penetration level of transmission-scale time-intermittent renewable generation resources increases, control of flexible resources will become important to mitigating the fluctuations due to these new renewable resources.  Flexible resources may include new or existing synchronous generators as well as new energy storage devices. Optimal placement and sizing of energy storage to minimize costs of integrating renewable resources is a difficult optimization problem.  Further, optimal planning procedures typically do not consider the effect of the time dependence of operations and may lead to unsatisfactory results.  Here, we use an optimal energy storage control algorithm to develop a heuristic procedure for energy storage placement and sizing.  We perform operational simulation under various time profiles of intermittent generation, loads and interchanges (artificially generated or from historical data) and accumulate statistics of the usage of storage at each node under the optimal dispatch. We develop a greedy heuristic based on the accumulated statistics to obtain a minimal set of nodes for storage placement. The quality of the heuristic is explored by comparing our results to the obvious heuristic of placing storage at the renewables for IEEE benchmarks and real-world network topologies.

}
\vspace{0.5cm}

\section{Introduction}
\label{sec:Intro}
Electrical grid planning has traditionally taken two different forms; operational planning and expansion or upgrade planning.  The first is concerned with the relatively short time horizon of day-ahead unit commitment or hour-ahead or five-minute economic dispatch.  The focus is on controlling assets that are already present within the system to serve loads at minimum cost while operating the system securely.  The second typically looks out many years or decades and is focused on optimal addition of new assets, with a focus on minimizing the cost of electricity over the long time horizon.  When a system consists entirely of controllable generation and well-forecasted loads, the network power flows do not deviate significantly or rapidly from well-predicted patterns.  In this case, expansion planning can be reasonably well separated from operational planning.  In the latter case, expansions may be optimized against only a handful of extreme configurations.

As the penetration of time-intermittent renewables increases, expansion and operational planning will necessarily become more coupled.  For an electrical grid with large spatial extent, renewable generation fluctuations at well-separated sites will be uncorrelated on short time scales\cite{Gibescu2009,08LBNL}, and the intermittency of this new non-controllable generation will cause the patterns of power flow to change on much faster time scales than before, and in unpredictable ways.  This new paradigm shift calls for accounting of multiple diverse configurations of uncertain resources in many operational as well as planning tasks. New equipment (e.g. combustion turbines or energy storage) and control systems may have to be installed to mitigate the network effects of renewable generation fluctuations to maintain generation-load balance. The optimal placement and sizing of the new equipment depends on how the rest of the network and its controllable components respond to the fluctuations of the renewable generation.  Overall, we desire to install a minimum of new equipment by placing it at network nodes where controlled power injection and/or consumption have a significant impact on the network congestion introduced by the renewable fluctuations.  From the outset, it is not clear which nodes provide the best controllability. Placing a minimum of new equipment is desirable since the investment and installation costs and costs associated with overcoming regulatory barriers. Thus, it makes sense to minimize the number of sites at which storage is placed for economic reasons.

\section{RELATED WORK}

Before discussing our initial approach at integrating operational planning and expansion planning, we summarize a few methods for mitigating the intermittency of renewable generation.  When renewable penetration is relatively low and the additional net-load fluctuations are comparable to existing load fluctuations, a power system may continue to operate ``as usual'' with primary and secondary regulation reserves\cite{99HK} being controlled via a combination of distributed local control, i.e. frequency droop, and centralized control, i.e. automatic generation control (AGC).  In this case, {\it planning} for renewables may simply entail increasing the level of reserves to guard against the largest expected fluctuation in {\it aggregate renewable output}.

As the penetration level grows, simply increasing the reserve levels will generally result in increased renewable integration costs\cite{Meibom2010} which are usually spread over the rate base.  Alternatively, operational planning can be improved by using more accurate techniques for renewable generation forecasting to better schedule the controllable generation (energy and reserves) to meet net load and operate reliably \cite{Meibom2010,Bouffard2008,BPA-wind-hirst}.   Simulations using rolling unit commitment \cite{Tuohy2007,Meibom2010}, where updated wind forecasts are used modify the unit commitment more frequently, have resulted in lower overall renewable integration costs.

Both unit commitment and economic dispatch seek to minimize the cost of electricity, however, they must also respect system constraints including generation/ramping limits, transmission line thermal limits, voltage limits, system stability constraints, and N-1 contingencies.  Previous works\cite{Meibom2010,Bouffard2008,Tuohy2007,BPA-wind-hirst} have generally looked at the effects of stochastic generation on the economics and adequacy of {\it aggregate} reserves while not considering such network constraints.  These constraints may be respected for a dispatch based on a {\it mean} renewable forecast. However, if the number of renewable generation sites and their contribution to the overall generation is significant, verifying the system security of all probable renewable fluctuations (and the response of the rest of the system) via enumeration is a  computationally intractable problem.

The approaches summarized above do not consider network constraints or the behavior of the system on time scales shorter than the time between economic dispatches (one hour in the case of \cite{Meibom2010}).  In particular, they do not model how fast changes in renewable generation and the compensating response of regulation reserves interact with network constraints.  In this manuscript, extending our initial study \cite{DvijothamChertkov}, we augment the approaches summarized above by focusing on the behavior of the electrical network at a finer time resolution and investigate how the control of energy storage affects its placement and sizing.

We presume that the unit commitment problem has been solved, and at the start of a time period, we perform time-varying (every 5 minutes) lookahead dispatch of controllable generation and storage based on an operational scenario (spatial and temporal profiles of wind generation, load and net interchange) while trying to minimize the storage capacity used (in terms of both energy and power) ---this gives us the minimum level of storage at each bus required for a particular operational scenario. We perform this optimization for several different scenarios (based on historical data, if available, or data generated using an appropriate statistical model). The statistics from simulated system operations are then coupled to the expansion planning process by developing a heuristic to guide the optimal placement and sizing of storage throughout the network---a result that cannot be achieved with the previous approaches described above.

A new approach, applying convex relaxations to traditional operations (like Optimal Power Flow (OPF)) including uncertain (wind) resources and storage was recently proposed in \cite{12BGTC}. The idea was to solve a version of the OPF problem with certain constraints relaxed (permitting potentially inadmissible solutions) so that the resulting problem is a convex optimization problem and can be solved to global optimality efficiently. Further, the authors provide conditions under which the solution to the relaxed problem satisfies all the constraints of the original problem, so that the relaxed problem can be used as a computationally efficient proxy. The approach was also extended to the storage placement problem in \cite{12BGTC,13GT}, which concluded, that placement of storage on, or close to, renewable sites is far from optimal. Although innovative and theoretically interesting, the convex relaxation approach of \cite{12BGTC,13GT} lacks scalability and was only illustrated on a very small 14 bus system. This is due to the high computational complexity of the semidefinite programming approach used in \cite{12BGTC}. Further, the authors in \cite{12BGTC} need to assume periodicity of renewable generation in order to solve the storage placement problem. In contrast, our work is the first resolving the storage placement problem over realistically sized networks.   We run our algorithm on a 2209-node model of the Boneville Power Administration (BPA), accounting for actual operational data and multiple (more than hundred) wind patterns.

As discussed earlier, there are several reasons to place energy storage at a small number of sites. However, choosing the optimal set of sites is a \emph{combinatorial} problem and cannot be solved by convex programming techniques.  In this paper, we develop a greedy heuristic that attempts to solve the storage placement problem directly. While we can no longer guarantee optimality of this algorithm, we demonstrate that our approach is robust and works across different network topologies leading to more economical placements that obvious alternatives.

The rest of this manuscript is organized as follows.  Section~\ref{sec:math} lays out the system model including the development of our heuristic algorithm for placement and sizing of energy storage.  Section~\ref{sec:exp} describes our numerical simulation results on a slightly modified version of RTS-96\cite{RTS96} and on a model of the BPA system. Section~\ref{sec:con} wraps up with some conclusions and directions for future work.

\section{Mathematical Formulation}
\label{sec:math}
We first set up notation in subsection
\ref{subsec:background}, formulate the look-ahead dispatch problem in Subsection \ref{subsec:OPF}, and finally describe our heuristic algorithm for storage placement in \ref{subsec:Sizing}.

\subsection{Background and Notation}
\label{subsec:background}

Let $\Gr=(\V,\E)$ denote the graph representing the topology of the power system, with $\V$ being set of buses and $\E$ being the set of transmission lines. Let $n$ denote the number of buses. At each bus, we can have three types of elements: Loads which consumer active power, Traditional Generators which generate active power and whose output can be controlled (within limits) and unconventional generators (renewables like wind) which generate active power, but whose output cannot be controlled.

For any $S\subset \V$ and any vector $v \in \mathcal{R}^{n}$, we denote $v_S=\{v_i:i \in S\}$. We will sometimes abuse this notation slightly to also denote the $n$ dimensional vector with zeros everywhere except in $S$. We denote by $\pv$ the vector of net-injections at each bus, and by $\pvs,\pvr,\pvl, \pvg$ the vector of injections at every node due to storage, renewables, loads and traditional generators, respectively. For any quantity $y$ that is a function of time, we denote by $y(t)$ its value at time $t$. In this paper, we will use integer-valued time $t=0,1,\ldots,T_f$ where $T_f$ is the time horizon of interest.

Let $\N(j)$ denotes the set of neighbors of node $j$ in the network. We define the graph Laplacian to be a $|V|\times |V|$ matrix with entries: \[\Lap_{ij}=-\frac{1}{\ind_{ij}},\Lap_{ii}=\sum_{j \in \N(i)}\frac{1}{\ind_{ij}}\]
 where $\ind_{ij}$ is the reactance on the transmission line between node $i$ and $j$. Then, the DC power flow equations are given by:
\[\fv_{ij}	=	\frac{\theta_i-\theta_j}{\ind_{ij}} ,\theta =  \Lap \pv\]
where $\theta$ denotes the voltage phase and $\fv_{ij}$ the active power flow between node $i$ and $j$.

We will consider placing energy storage at nodes in the network. We denote the by $\sv$ the vector of energy stored at each node in the network. The energy capacity of storage (maximum energy that can be stored) is denoted $\svlim$ and the maximum power that can be withdrawn from or supplied to the energy storage units $\pvslim$. We denote by  $\pvlim$ the maximum power output of traditional generators and $\pvglim$ the corresponding limit on the ramping limits. $\flim_{ij}$ the limit on the flow on the line between $i,j$.

\subsection{Lookahead Dispatch of Generation and Storage}
\label{subsec:OPF}

In the presence of energy storage, the Optimal Power Flow(OPF)-based dispatch problem gets coupled over time (since energy stored at some time can be used later). Our approach to sizing and placing energy storage relies on operational simulations of the system under realistic load and renewable generation profiles. The operational simulation is formulated as a lookahead-dispatch problem: This is very similar to what the system operator would do to dispatch energy storage given a forecast of renewable generation and load. However, since we are interested in sizing and placement of energy storage, we additionally optimize over the energy capacity $\svlim$ and power capacity $\pvslim$ of the energy storage needed to ameliorate the fluctuations in renewables and loads.

\begin{align}
 &\min_{\pvs(t),\pvg(t)} \underbrace{\sum_{t=0}^{T_f} \tran{\costg}\pvg(t)}_{\text{Generation Costs}}
 + \underbrace{\tran{\costs}\svlim+\tran{\costsp}\pvslim}_{\text{Storage Investment Costs}}
  &\nonumber \\
 &\textrm{subject to} & \nonumber \\
 & 0 \leq \pvg(t) \leq \pvlim_g \quad \text{(Generation Capacities)}  \label{eq:OptC1}\\
 &\pv(t) = \pvg(t)+\pvr(t)+\pvl(t)+\pvr(t) \quad \text{(Net Injection)} \label{eq:OptC2} \\
 &\Lap\theta(t) = \pv(t) \quad \text{(DC Power Flow)}  \label{eq:OptC3}\\
 &\left|\fv_{ij}(t) = \frac{\theta_i(t)-\theta_j(t)}{\ind_{ij}}\right|\leq \flim_{ij} \quad \text{(Flow Limits)} \label{eq:OptC4} \\
 &\left|\pvg(t+1)-\pvg(t)\right|\leq \pvglim \quad \text{(Generation Ramping Limits)}  \\
 & 0 \leq \sv(t) \leq \svlim \quad \text{(Energy Capacity of Storage)}  \label{eq:OptC5}\\
 & 0 \leq \pvs(t) \leq \pvslim \quad \text{(Power Capacity of Storage)}  \label{eq:OptC6}\\
  &\sv(t) = \sv(0)-\sum_{\tau=0}^{t-1} \pvs(t)\Delta  \quad \text{(Energy Conservation)} \label{eq:OptC7}\\
  &\tran{\One}\sv(T_f)=\tran{\One}\sv(0) \quad \text{($0$ Net Energy Supply)}  \label{eq:Opt}
\end{align}

The objective models operational costs of generation (fuel etc.) and \emph{amortized} investment costs of placing energy storage in the grid. The constraints \eqref{eq:OptC1},\eqref{eq:OptC2},\eqref{eq:OptC3} and \eqref{eq:OptC4} are standard constraints appearing in a DCOPF formulation. The fifth constraint \eqref{eq:OptC5} is relevant in scenarios where wind generation undergoes a ramp event (sudden drop or increase) and traditional generators need to increase or decrease their output at rates close to their ramping limits. The constraints \eqref{eq:OptC5}, \eqref{eq:OptC6}, \eqref{eq:OptC7} are standard constraints for storage. The final constraint \eqref{eq:Opt} models the fact that we want to use energy storage as a hedge over time - to store energy when too much power is being produced in the grid and supply it at a later time. Thus, over the horizon of interest, we do not want a net energy supply to/from the energy storage. This optimization problem is a Linear Program (LP) (like a standard DCOPF) and can be solved using off-the-shelf linear programming packages. We use the gurobi package in our work here \cite{gurobi}.

\subsection{Modeling Assumptions} \label{sec:DCvsAC}
Since this is a preliminary study meant to illustrate the value of coupling planning and operations, we made a number of simplifying assumptions that may not hold for a real power system. The first one is to use the DC Power Flow equations rather than the full nonlinear AC equations. The second one is to assume that dispatch is based on perfect forecasts of wind and loads over a 2-hour period. We outline the justifications for these assumptions in this section.

\subsubsection{DC vs AC OPF}
The DC power flow equations are an approximation to the nonlinear AC power flow equations. They are  frequently used in the context of power markets although system operators would use the nonlinear ACOPF to perform actual dispatch of generators in a grid. In general, there can be significant discrepancies between DC and AC power flow results that make the DC solution unacceptable in an operational setting. In this paper, however, we stick with the DCOPF formulation. There are multiple reasons for this:

\begin{itemize}

\item[1] Since our interest in this work was to concentrate on the novel aspect of integrating planning and operational studies, we were not interested in building a nonlinear ACOPF solver. Freely available solvers like MATPOWER \cite{zimmerman2005matlab} do not generalize to the lookahead dispatch setting, that is, they are unable to deal with the time-coupling introduced by storage. However the storage placement algorithm (Algorithm \ref{alg:OptPlace}) we develop in this paper can be used with any OPF solver. In particular, a more complete commercial-grade ACOPF solver should work better. Additionally, the extra computational burden of the ACOPF is not an issue here since we are performing offline planning studies which does not impose strict real-time requirements on the computation time (we could allow the algorithm to run for days if required).

\item[2] We are mostly concerned with \emph{long-term planning} and use operational information to inform the planning process. Hence, we are only interested in the accuracy of the OPF to the extent that it captures all possible patterns of flows observed in typical operational scenarios. In numerical studies we performed, we observed that the DCOPF suffices for this purpose, at least for this preliminary study meant to illustrate the value of coupling planning and operations.

\item[3] In general, the DCOPF becomes less accurate as the system gets under more stress. While we consider systems with high penetrations of renewable energy, we do not aim to deal with critical scenarios where the grid is under stress (close to voltage/ frequency instability). The challenge of high renewable penetration (which we aim to handle here) is that of non-predictable patterns of power flows. Thus, we are looking at the system under stable operating conditions, but with fluctuating patterns of power flows. When the grid is under stress, we assume that appropriate emergency control actions will be taken to protect the system. We do not aim to use energy storage to perform emergency control actions.
\end{itemize}

\subsubsection{Perfect Forecasts}
Note that in our DCOPF formulation $\pvr(t),\pvl(t)$ are assumed to be known functions of time. This is like performing lookahead dispatch with perfect forecasts. Although this differs from a real operational scenario (imperfect forecasts), we believe that the discrepancy will not break our analysis here for the following reasons:

\begin{itemize}
\item We consider time horizons of about $T_f=2$ hours. Over such a time-scale, loads are well-predictable for sure, although wind may not be. However, we use operational simulations to develop a heuristic for placement of energy storage: Hence changes due to forecast errors, while important in an operational context, are less important from the context of deciding placement of energy storage.
\item Several system operators today perform periodic redispatch of the grid resources (generation/storage) at fairly short intervals of time (5-15 mins) and hence can easily adapt to and cope-with forecast errors.
\end{itemize}

Further, we note that our heuristic for storage placement is \emph{independent} of the specific dispatch scheme (OPF) used. Thus, we can perform a robust or chance-constrained version of DCOPF \cite{12BCH} which would allow us to incorporate the effect of forecast uncertainty into the dispatch, and hence into the storage placement decision.

\subsection{Optimal Sizing and Placement of Storage}\label{subsec:Sizing}

We seek to develop heuristics to decide how to place storage and size its energy and power capacity. However, we must first define some metrics to evaluate a given storage placement.  Let $\So$ denote the set of nodes with non-zero storage. For a given scenario $\delta_i$ (renewable/load profiles) and $\So$, the energy and power capacities resulting from the optimization \eqref{eq:Opt} are $\svlim_i$ and $\pvslim_i$. We define the energy in the renewable fluctuations to be $\sr(t)=\sum_{\tau=0}^{t-1} \pvr(\tau)\Delta$, i.e. $\sr(t)$ is the energy stored in a (hypothetical) battery that is connected directly to a renewable node and eliminates all fluctuations about the mean renewable generation. Then, plausible metrics can be defined according to the following criteria:

\emph{Normalized Power Capacity}: This quantifies the total power capacity of the storage relative to the sum of maximal power fluctuations over the renewables:
\[\frac{\sum_{j \in \So} \max_t\abs{\popt_j(t)}}{\sum_{i} (\max_t \pvr_i(t)-\min_t\pvr_i(t)) }\]
\emph{Normalized Energy Capacity}: This quantifies the total energy capacity of the storage relative to the sum of maximal energy fluctuations over the renewables:
\[\frac{\sum_{j \in \So} (\max_t \sopt_j(t)-\min_t \sopt_j(t))}{\sum_{i} (\max_t \sr_i(t)-\min_t\sr_i(t)) }\]
\emph{Overall Performance}: We denote a weighted combination of the above metrics by $\mathrm{perf}(\So)$. In this study, we choose this to be the total normalized energy capacity plus a fixed cost for each site at which storage needs to be placed.
\\
\emph{Renewable Penetration}: The fraction of load served by renewables over the time horizon $T$.\\

 The high-level pseudocode given in Algorithm~\ref{alg:OptPlace}.  The algorithm is a greedy pruning heuristic that starts with $\So=\V$, i.e.  storage at all nodes, and seeks to shrink $\So$ while improving performance at least by some minimum amount $\epsilon$ at each iteration.  Then, the same procedure is repeated, each time shrinking the target number of nodes as long as the performance metric is improving.
 Note that this repetition is required (and critical) because the dispatch based on restricted storage would be different, since there are a smaller set of controllable resources.
\begin{algorithm}
\caption{Greedy Heuristic for Optimal Placement}  \label{alg:OptPlace}
\begin{algorithmic}
 \STATE Input: Collection of Scenarios $\{\delta_k\}$, Threshold $\epsilon$
  \STATE $\So \gets \{1,2,\ldots,n\}$.
 \REPEAT{}
 \FOR{$k = 1 \to N$}
  \STATE Solve \eqref{eq:Opt} for scenario $\delta_k$ to get $\svlim_k,\pvslim_k$
 \ENDFOR
  \STATE $\svlim \gets \max_{k} \svlim_k$
  \STATE $\pvslim \gets \max_{k} \pvslim_k$
  \STATE $\gamma \gets \max\{\gamma:\{\mathrm{perf}(\{i \in \So:\svlim_i\geq\gamma\max(\svlim)\})<\mathrm{perf}(\So)-\epsilon\}\}.$
  \STATE $\So \gets \{i \in \So:\svlim_i\geq\gamma\max(\svlim)\}.$
 \UNTIL{$1-\gamma\leq\epsilon'$}
\end{algorithmic}

\end{algorithm}

\subsection{Justification for Greedy Algorithm}

The choice of the greedy algorithm is motivated by the theory of submodular function maximization \cite{krause2012submodular}. Submodular functions are functions with diminishing marginal returns. Mathematically, if one had a function $F$ defined on subsets $A$ of $S=\{1,2,\ldots,m\}$ that satisfied:
\[F(A \cup \{i\}) - F(A) \leq F(B \cup \{i\}) - F(B), B \subset A \subset S, i \not\in A.\]

In our context, this simply means: The additional performance gain obtained by adding storage at a new node when there is already storage at a large number of nodes is smaller than the performance gain obtained by adding to storage when there is storage at only a few nodes. Although we have not been able to prove that this property holds for the storage placement problem, it definitely makes intuitive sense - at some point one would expect to observe diminishing returns for additional placement. It can be shown that for a submodular objective function, the greedy algorithm achieves an objective that is within $1-\frac{1}{e}$ of the optimal solution \cite{krause2012submodular}. This motivated us to consider a greedy algorithm to solve the problem of storage placement.

We have some preliminary results (not included in this paper) regarding the submodular property for certain simplified versions of the objective presented here and hope to pursue this line of investigation further in future work.

\section{Simulations}
\label{sec:exp}
\label{sec:Results}
\subsection{RTS-96+Synthetic Wind Data}
\begin{figure}[htb]
\centering
\subfigure[Our modified version of RTS-96.  The added renewables are blue, loads are yellow and controllable generators are green]{
 \includegraphics[width=.3\textwidth]{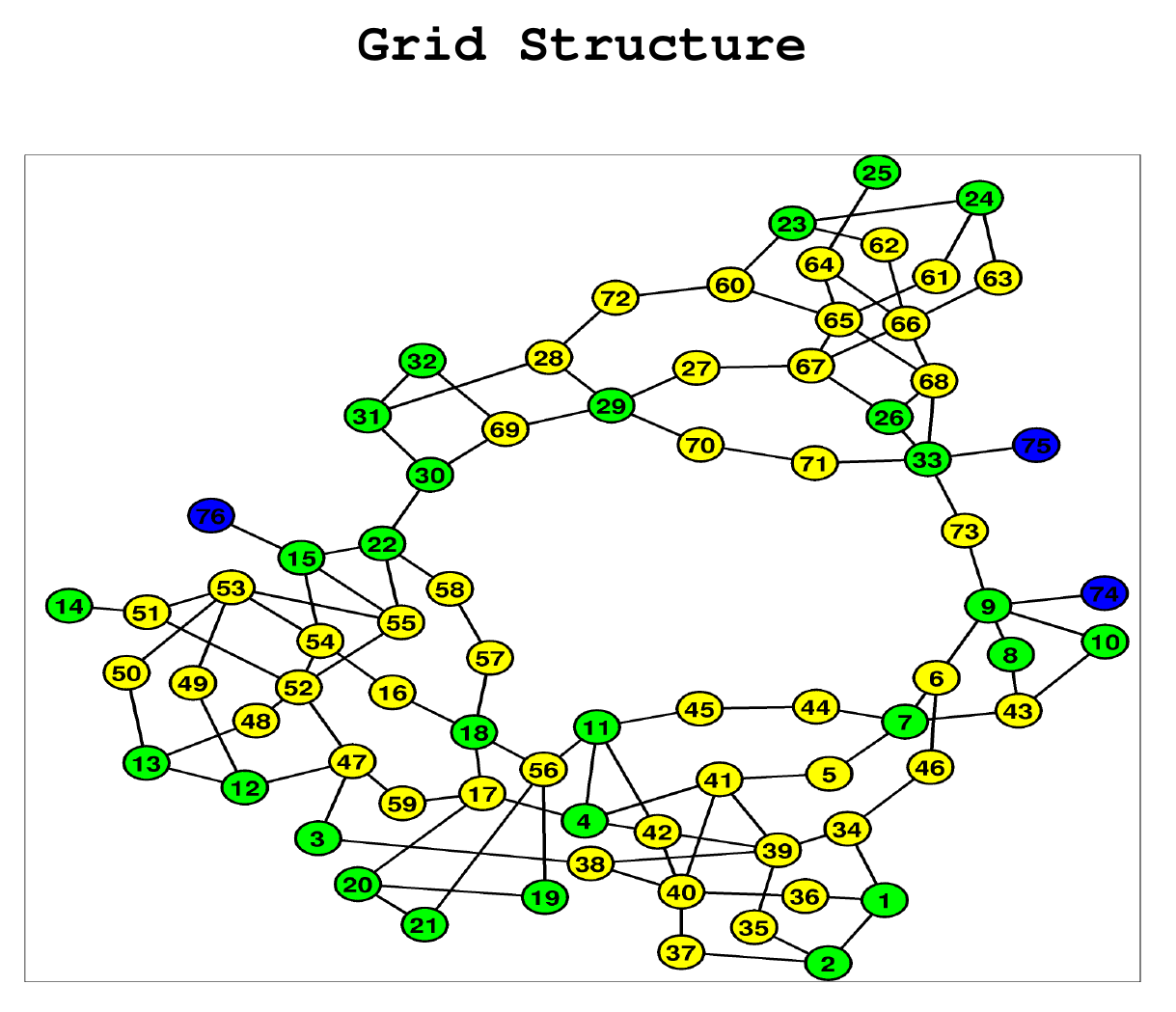}
 \label{fig:Grid}}
\subfigure[Sets of nodes identified by our heuristic: The minimal set (two nodes) is shown in red, additional 8 nodes, of the 10 node set, are show in green, all other nodes are shown in blue.]{
 \includegraphics[width=.3\textwidth]{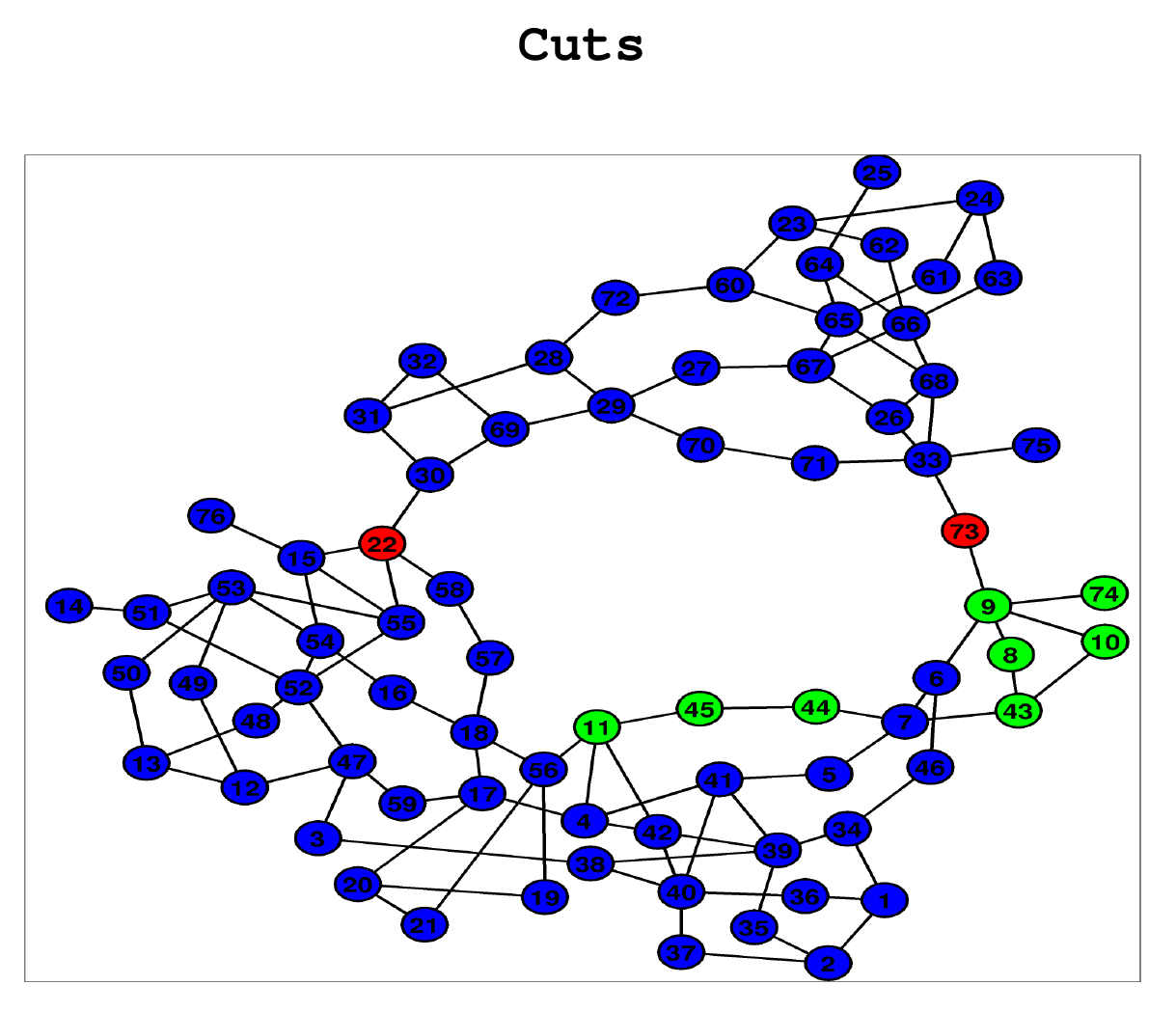}
\label{fig:Cuts}}
\subfigure[Storage Capacity Histograms: Red lines mark thresholds used for the reduction in the storage node set]{
 \includegraphics[width=0.3\textwidth]{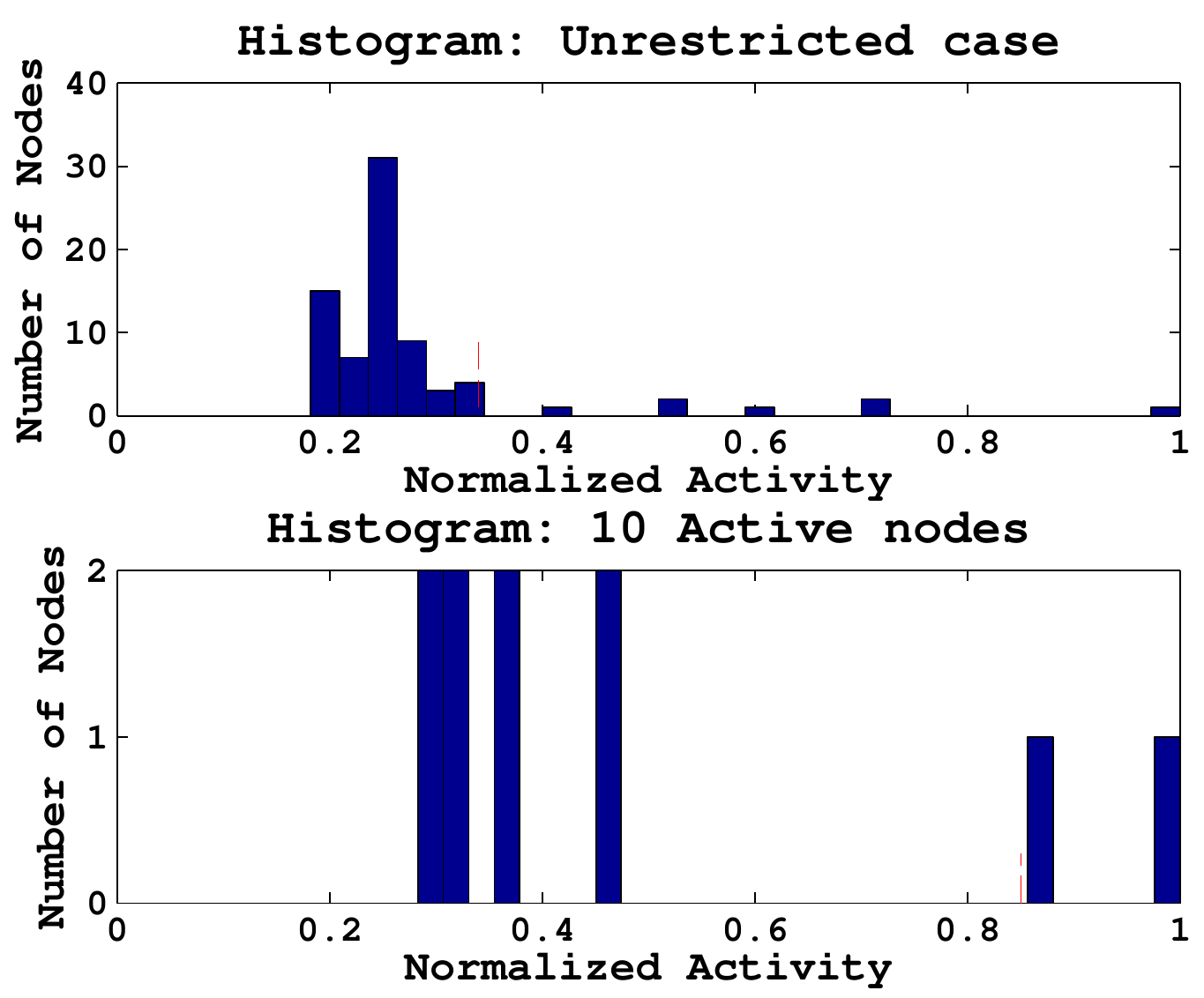}
\label{fig:Hist}}
\end{figure}

\begin{figure}[htb]
\centering
\subfigure[The normalized energy capacity of storage in the entire network vs the penetration of renewable generation]{
 \includegraphics[width=0.35\textwidth]{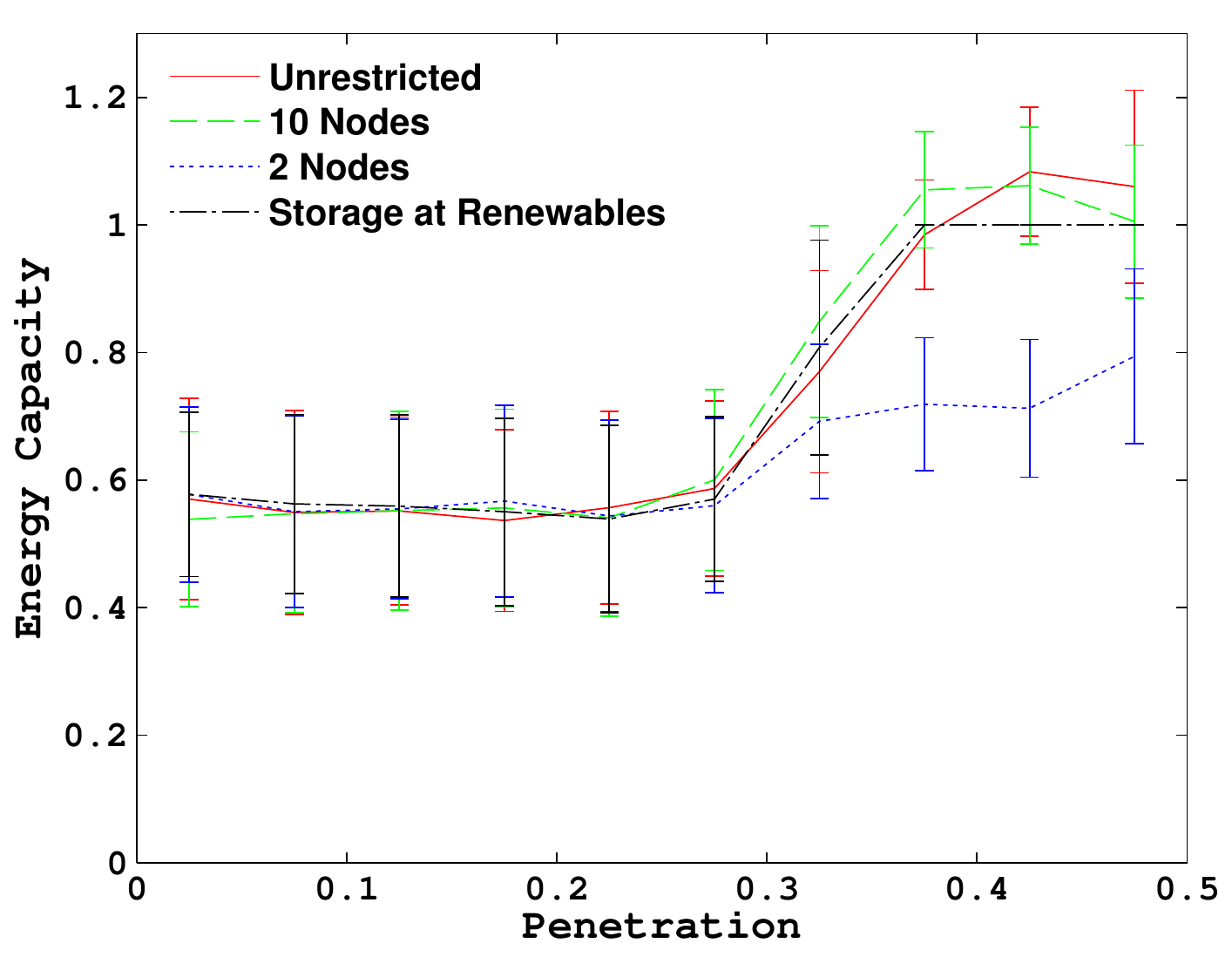}
\label{fig:EnergyCapP}}
\subfigure[The normalized power capacity of storage in the entire network vs the penetration of renewable generation]{
 \includegraphics[width=0.35\textwidth]{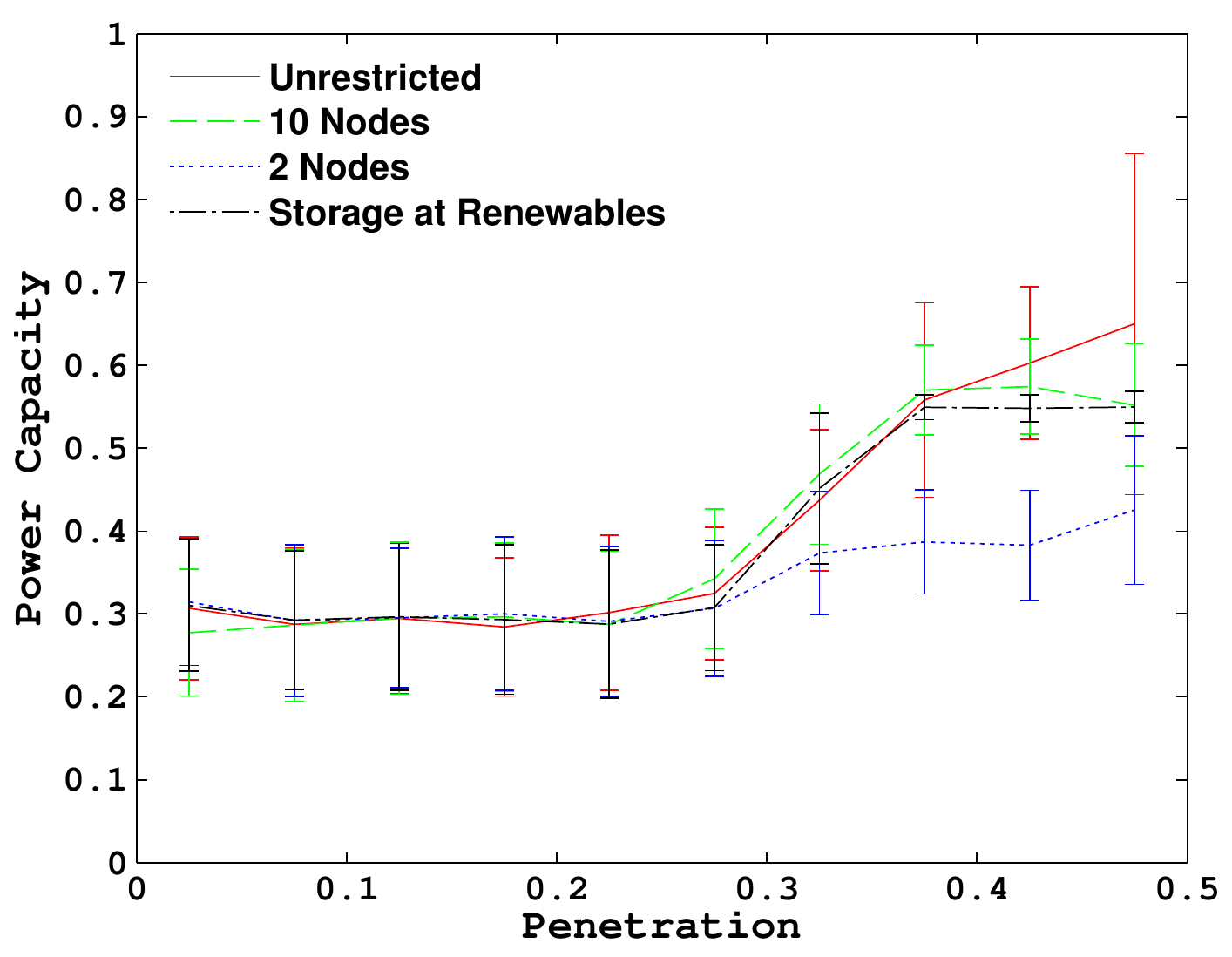}
\label{fig:PowerCapP}}
\end{figure}

We tested our optimal control and heuristic for storage placement and sizing on a modified versions of RTS-96\cite{RTS96}.  The grid is shown in Fig.~\ref{fig:Grid}.  Our modification includes the addition of three renewable generation nodes shown Fig.~\ref{fig:Grid} in blue.  The capacities of the new lines connecting the renewables to their immediate neighbors are set higher than the capacity of the added renewable generation, otherwise, these lines would be overloaded in nearly every trial.

In each iteration of algorithm \ref{alg:OptPlace}), we generate $N=2000$ time series profiles for the renewables. These are chosen so that we can control the penetration of wind in the system and study the effect of penetration of intermittent sources on storage sizing and placement.

In the first iteration, storage is available at all nodes in the network.  The histogram of storage capacities is plotted in Fig.~\ref{fig:Hist}. We then shrink the set of nodes having storage until the performance metric $\mathrm{perf}({\bf S}) $ (defined in Section \ref{subsec:Sizing}) fails to improve significantly (by more than $\epsilon$). For this example, we were able to shrink down to $10$ nodes in the first iteration. Using these 10 nodes, we rerun the optimal control algorithm and again accumulate statistics of the storage activity (plotted in figure \ref{fig:Hist}). Based on the updated statistics, we can again shrink the set of storage nodes down to $2$ nodes and this is the final output of the algorithm (we cannot shrink any further without performance degradation). The optimally chosen sets of 10 and then 2 nodes are shown in Fig.~\ref{fig:Cuts}.

The method for generating the renewable profiles is described in details in \cite{DvijothamChertkov}. The evaluation metrics defined in Section \ref{subsec:Sizing} are shown as functions of penetration in Figs.~\ref{fig:EnergyCapP},\ref{fig:PowerCapP}, with storage at all the nodes and sets of shrunken nodes discovered by the Algorithm \ref{alg:OptPlace}.

\subsection{BPA System with Historical Wind Data}

We also apply our algorithm to real data from the BPA network covering Washington and Oregon. By overlaying the grid on the US map, we were able to locate the major wind farms and inter-ties (to California) in the system. Loads were divided roughly in proportion to population densities. Mapping this onto data published on the BPA web-site \cite{BPA_wind}, we were able to create realistic wind, load and interchange profiles. We considered data from 100 different wind configurations during 2012 (each of length about 2 hours, spread uniformly throughout the year). We also ensured that we pick particularly challenging operating conditions, for example, periods with high ramping conditions in wind generation, i.e. these pushing the storage dispatch to its limits, and thus to enable sizing storage so as to be prepared for the worst contingencies.  

\begin{figure}
\centering
\subfigure[BPA - Algorithm Iteration 1]{
 \includegraphics[width=.3\textwidth]{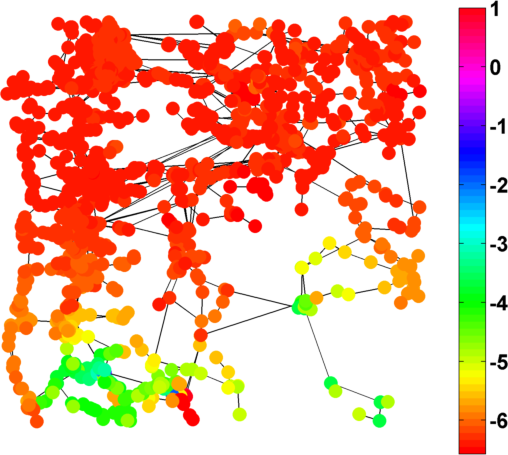}
 \label{fig:BPA1}}
   \subfigure[BPA - Algorithm Iteration 2]{
 \includegraphics[width=.3\textwidth]{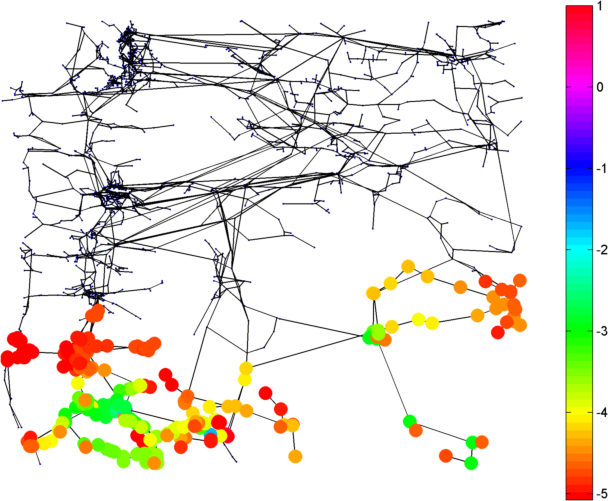}
 \label{fig:BPA2}}
 \subfigure[BPA - Algorithm Iteration 3]{
 \includegraphics[width=.3\textwidth]{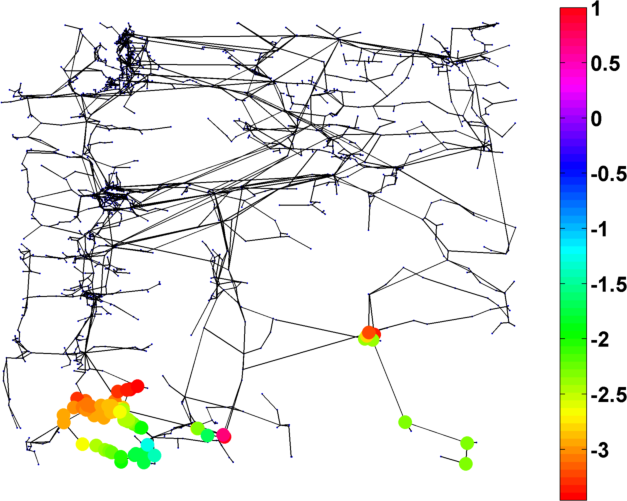}
 \label{fig:BPA3}}  \subfigure[BPA - Algorithm Iteration 4]{
 \includegraphics[width=.3\textwidth]{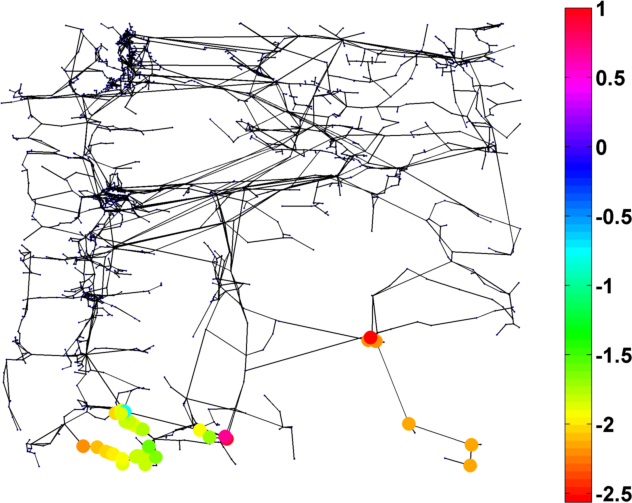}
 \label{fig:BPA4}}
\caption{Iterations of our Algorithm on the BPA System. Red Corresponds to Low Storage Capacity and Purple to High}\label{fig:AlgBPA}

\end{figure}

\begin{figure}
 \includegraphics[width=.4\textwidth]{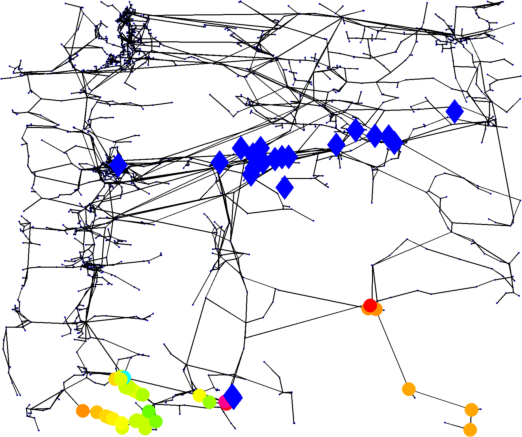}

\caption{Storage Placement (colored circles) relative to Wind Farms/Interties (shown as blue diamonds).
} \label{fig:BPAWind}

\end{figure}

We plot the iterations of our algorithm on the BPA system in Fig.~\ref{fig:AlgBPA}. The nodes at which storage is present are colored---red marking the nodes with least storage capacity and purple marking the nodes with the highest storage capacity - The capacities are color coded in a log-scale:
\[\log\left(\frac{\text{Storage Capacity at a node}}{\text{Maximum Storage Capacity over all nodes}}\right)\]
so as to improve visual discriminability. Our sequential algorithm is able to discover a relatively small subset of 37 nodes at which to place storage. Reducing this number any further leads to a significant increase in the overall storage capacity required.

\begin{figure}
   \begin{center}

    \includegraphics[width=.34\textwidth]{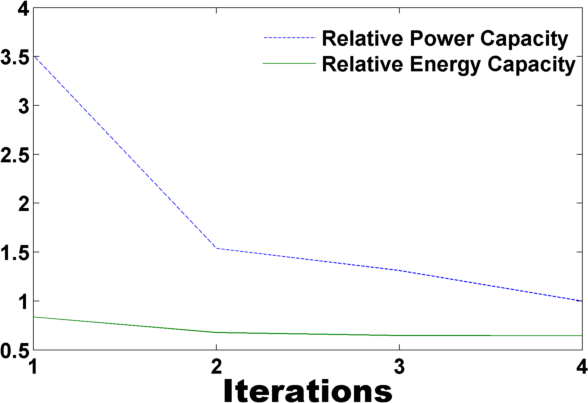}

    \end{center}

           \caption{Total Energy and Power Capacity Relative to Placement at Renewables and Interties}\label{fig:Iter}

   \end{figure}

In Fig.~\ref{fig:BPAWind}, we plot the locations of the storage nodes relative to the locations of the wind farms and inter-ties.   We note that our algorithm does not place the storage near either the wind sites or the interties  We also compared our strategy to placing storage directly at the wind farms or inter-ties (which are the ``sources" and "sinks" which contribute most to fluctuations in the generation/load). The overall storage capacity required by this naive approach is twice the storage capacity required for the placement discovered by our algorithm.  In Fig.~\ref{fig:Iter},we plot the total energy and power capacity of the storage placements discovered by our algorithm relative to the naive strategy of placing storage directly at the renewables and interties. 

This result shows that the storage placement discovered by the algorithm, although intuitive, is non-trivial since the many of the nodes picked are not precisely at the renewables or interties, but rather at critical nearby nodes which are critical for controlling the power flows in that region of the network. 

\subsection{Computation Times}

For the BPA system, the entire algorithm took about 10 minutes on a i7 2.9 GHz CPU to produce the optimal placement of storage. The optimal dispatch for a particular scenario takes about 5 seconds.

\section{Discussion}

We have presented an efficient and effective heuristic for sizing and placing energy storage in a transmission network. Our essential insight in this paper was to couple operational simulations with planning, and use statistics from operational simulations to inform the planning procedure. For any realistic engineered network, operational simulations will contain valuable information about the various flow patterns, congestion, ramping restrictions etc. in the network and provide an effective heuristic for making planning decisions - we have observed this in the above simulations as well. With an unoptimized matlab implementation, our approach takes about 10 minutes to discover an effective storage placement for the BPA system. For an offline planning problem, this is perfectly acceptable.

An alternate approach would be to formulate this directly as a mixed integer linear program: Choose a small number of sites to place energy storage so as to minimize investment and operational costs over a large set of possible scenarios. However, this approach fails to take advantage of the above observation, and quickly becomes computationally infeasible for realistically sized networks.

We use the DCOPF approximation in our work, but as mentioned in section \ref{sec:DCvsAC}, the approach can be easily used with an ACOPF solver. 

For both the BPA and RTS-96 network, by using the greedy pruning algorithm,  we are able to reduce the number of energy storage sites to a very small number compared to the total number of nodes in the network.  In both cases,  the storage is placed far from the renewable generation.   Instead,  the storage appears to be placed at a few critical nodes suggesting that the storage is being used not only to buffer fluctuations, but also to assist with controlling flows in the rest of the network. For the BPA system, these may seem geographically close to the renewables or interties. However, the precise placement of storage is non-trivial and the discovered placement uses particular nodes that offer a large degree of controllability on the power flow patterns in that region of the network. Thus, in effect,  our algorithm is designing the grid control system by finding the nodes with the highest controllability over the network congestion.

This conclusion is supported by the plot of iteration-by-iteration storage energy and power capacity in Fig. 3.  The energy capacity of the storage is not dramatically reduced by during the pruning.  Instead,  storage capacity that was dispersed throughout the network is concentrated at fewer nodes resulting in larger but sparser storage installations.  However,  the storage power capacity drops significantly.   This seems to indicate that the wind fluctuations require a certain amount of energy capacity for buffering on a network wide basis.   However, better placement of that energy capacity enables it  to be used just as effectively with a much smaller power capacity.

\section{Conclusions and Future Work}
\label{sec:con}

Somewhat unexpectedly, our algorithm chooses to place storage at nodes at critical junctions between major subcomponents of the network rather than at the sites of renewable generation. We conjecture that these nodes provide for enhanced controllability because, in addition to simply buffering the fluctuations of the renewables, controlled power injections at these nodes can modify overall network flows and direct fluctuating power flows to regions that are better positioned to mitigate them.

There is much follow on work needed to expand the concept presented in this manuscript and to verify some of its conjectures. We also plan to extend this approach to allow for stochastic, robust and/or chance-constrained optimization, as in \cite{12BCH}, to provide for a better representation and more accurate modeling of the wind uncertainty. Finally, we performed lookahead dispatch assuming perfect information about loads and renewable generation. Although this assumption is reasonable (for reasons described in section~\ref{subsec:OPF}), a more thorough study is required to determine the exact effect of forecast errors, particularly on storage sizing (we would expect the placement to be robust to reasonable forecast errors). This would require modeling standard generation response mechanisms (primary control and AGC) which modify generator outputs in response to changes in renewable generation and loads. These mechanisms are well studied for generation, but need to be extended for storage systems as well. We plan to build on recent work in this direction  \cite{12DBC} for this.



{\small
\bibliographystyle{IEEEtran}
\bibliography{SmartGrid}
}

\end{document}